\documentclass[12pt]{amsart}
\usepackage{amscd, amssymb}

\newtheorem{lm}{Lemma}
\newtheorem{tm}{Theorem}

\newcommand{\proj}{\mathbf P}
\newcommand{\grass}{\mathbf G}
\newcommand{\barr}{\overline}
\newcommand{\rarr}{\rightarrow}
\newcommand{\oh}{{\mathcal{O}}}
\newcommand{\com}{\mathbb{C}}
\newcommand{\Q}{\mathbb{Q}}
\newcommand{\Z}{\mathbb{Z}}
\newcommand{\G}{\mathbf{G}}
\newcommand{\Or}{\mathbf{O}}
\newcommand{\SO}{\mathbf{SO}}

\newcommand{\PGL}{\mathbf{PGL}}
\newcommand{\cal}{\mathcal}
\newcommand{\eqq}{\stackrel{\sim}{=}}

\newcommand{\sumo}{\oplus}

\newcommand{\bpf}{\noindent {\em Proof.} }
\newcommand{\epf}{\qed \vspace{+10pt}}

\begin{document}
\begin{center}{\bf The Equivariant Chow Ring of
$\SO(4)$}
\end{center}
\begin{center}{Rahul Pandharipande$^1$}
\end{center}
\begin{center} {4 July 1996}
\end{center}

\pagestyle{plain}
\footnotetext[1]{Research 
partially supported by an NSF Post-Doctoral Fellowship.}
\setcounter{section}{-1}
\section{\bf{Introduction}}
Let $\G$ be a reductive algebraic group. 
The algebraic analogue of $E\G$ is attained by
approximation. Let $V$ be a $\com$-vector space. Let
$\G\times V \rarr V$ be an
algebraic representation of $\G$. 
Let $W\subset V$ be a $\G$-invariant open set satisfying:
\begin{enumerate}
\item[(i)] The complement of $W$ in $V$ is of codimension greater than 
$q$.
\item[(ii)] $\G$ acts freely on $W$.
\item[(iii)] There exists a geometric quotient $W\rarr W/\G$.
\end{enumerate}
$W$ is an approximation of $E\G$ up to codimension $q$.
Let $e=dim( W/{\G})$.
The equivariant Chow groups of $\G$ (acting on a point) 
are defined by:
\begin{equation}
\label{defff} 
A^{\G}_{-j}(\text{point})= A_{e-j}(W/\G)
\end{equation}
for $0\leq j \leq q.$
An argument is required to check the
Chow groups are well-defined (see [EG1]).
The basic properties of equivariant
Chow groups are  
established in [EG1]. In particular, 
there is a natural intersection
ring structure on $A_i^{\G}(\text{point})$.
For notational convenience, a superscript will
denote the Chow group codimension:
 $$A^{\G}_{-j}(\text{point}) = A^j_\G(\text{point}).$$
Equation (\ref{defff}) becomes:
$$\forall\  0\leq j \leq q, \ \
A_{\G}^{j}(\text{point})= A^j(W/\G).$$
$W/\G$ is an approximation of $B \G$.
$A_{\G}^*(\text{point})$ is called the (equivariant) Chow ring of
$\G$. 
$A_{\G}^*(\text{point})$ is naturally isomorphic to the
ring of algebraic characteristic classes of (\'etale locally
trivial) principal $\G$-bundles.
The equivariant Chow ring of $\G$ was
first defined by B. Totaro in [T].

Consider now the
orthogonal and special orthogonal algebraic groups
(over $\com$).
The Chow ring of $\Or (n)$ is    
generated by the Chern classes
of the standard representation. The odd classes
are 2-torsion:
$$A^*_{\Or (n)}(\text{point})= 
\Z[c_1, \ldots, c_n]/(2c_1, 2c_3, 2c_5, \ldots).$$
The Chow ring of $\SO (n=2k+1)$
is also generated by the Chern classes
of the standard representation. The odd classes
are 2-torsion and $c_1=0$:
$$A^*_{\SO (n=2k+1)}(\text{point})=
 \Z[c_1, \ldots, c_n]/(c_1, 2c_3, 2c_5, \ldots).$$
$A^*_{\Or (n)}(\text{point})$ was first computed
by B. Totaro. Algebraic computations of
$A^*_{\Or (n)}(\text{point})$ and $A^*_{\SO (2k+1)}(\text{point})$
can be found in [P2]. The Chow ring of
$\SO(n)$ has been computed with $\Q$-coefficients
in [EG2]. The integral Chow ring of
$\SO(n=2k)$ is not known in general. 

Since 
$B \Or (n)$ is approximated by the set of
non-degenerate quadratic forms in $Sym^2 S^*$ (where
$S \rarr \grass(n, \infty)$ is the
tautological sub-bundle), the Chow ring
can be analyzed by degeneracy calculations ([P2]). 
Algebraic
$B \SO (n)$ double covers $B \Or(n)$.
If $n=2k+1$, there
is a product decomposition $\Or(n)\eqq \Z/2\Z \times \SO(n)$.
As a result,
the double cover geometry for $n=2k+1$ is tractable and
a computation of $A_{\SO(2k+1)}^*(\text{point})$ can be made. 
In case $n=2k$, the double cover geometry is more
complicated.

Since
$\SO(2)\eqq \com^*$, the first non-trivial even case is
$\SO(4)$.
In this paper, 
the ring $A^*_{\SO(4)}(\text{point})$ is determined.
Let $c_1,c_2,c_3,c_4$ be Chern classes of the
standard representation of $\SO(4)$. Let $F$ be one of the
two distinct irreducible $3$-dimensional representations of $\SO(4)$, and
let $f_2$ be the second Chern class $F$.
\begin{tm}
$A_{\SO(4)}^*(\text{\em point})$ is generated by the Chern classes
$c_1, c_2, c_3,$ $c_4,$ and $f_2$. Define $x\in A_{\SO(4)}^{2}$ by
$x=c_2-f_2$.
 $$A_{\SO(4)}^*(\text{\em point})= \Z[c_1,c_2,c_3,c_4,x]/
(c_1, 2c_3, xc_3, x^2-4c_4)$$
\end{tm}
\noindent Let $\tilde{F}$ be the other irreducible $3$-dimensional
representation of $\SO(4)$. Let $\tilde{f}_2$ be the
second Chern class of $\tilde{F}$. Since (see [FH])
$$F \oplus \tilde{F} \eqq \wedge^2 V,$$
the relation $\tilde{f}_2= 2c_2-f_2$ is obtained.
Hence, $c_2-\tilde{f}_2= -x$. The presentation in
Theorem 1 does not depend upon the choice of $3$-dimensional
representation.

Thanks are due to D. Edidin, W. Fulton, W. Graham, and B. Totaro for
conversations about $B \SO(n)$. The $\SO (4)$ calculation
presented here is similar in spirit to the $\SO(2k+1)$
calculations of [P2]. In [P1] and [P2], equivariant Chow
rings are used to compute ordinary Chow rings of certain
moduli spaces of maps and Hilbert schemes of rational curves.

\section{\bf{Ruled Quadric Surfaces}}
\label{rrr}
Let $V\eqq \com^4$ be equipped with
a non-degenerate quadratic form $Q$. A {\em ruled quadric
surface} in $\proj(V)$ is a pair $(X,r)$ where  $X\subset \proj(V)$
is a nonsingular quadric surface and $r$ is a choice of ruling.
Let
$\cal{X}\subset \proj(Sym^2 V^*)$ be the parameter space
of nonsingular quadrics.
The  parameter space of ruled quadrics, 
$\cal{X}_{ruled}$, is an \'etale double
cover of $\cal{X}$ via the natural map:
$\cal{X}_{ruled} \rarr \cal{X}$.

There are natural maps $\SO(V) \rarr \proj\SO(V)\subset \PGL(V)$
and
$\Or(V) \rarr \proj\Or(V) \subset \PGL(V)$. 
Let $\SO(V)$ and $\Or(V)$ act on $\PGL(V)$ on the
right via these maps.
There exist geometric
quotients (see [P2]):
$$\PGL(V)/\SO(V) \rarr \PGL(V)/\Or(V).$$
Consider the quadric surface $(Q)\subset \proj(V)$
obtained from the quadratic form.
The standard left action $\PGL(V)\times \proj(V) \rarr
\proj(V)$ yields a  transitive $\PGL(V)$-action on the
space of nonsingular quadric surfaces. The
stabilizer of $(Q)$ for this action is exactly 
$\proj\Or(V) \subset \PGL(V)$. Hence, there is a
canonical isomorphism
$$\PGL(V)/\Or(V) \eqq \cal{X}.$$ 
For the entire paper, fix a ruling $r$ of $(Q)$. 
Since $\PGL(V)$ acts transitively on the
space of ruled quadrics and the stabilizer of
$((Q),r)$ is exactly $\proj\SO(V) \subset \PGL(V)$,
there is a canonical isomorphism determined by $((Q),r)$: 
$$\PGL(V)/\SO(V) \eqq \cal{X}_{ruled}.$$

There is a canonical Pl\"ucker embedding
$\grass(2, V) \hookrightarrow \proj(\wedge^2 V)$.
Let $Z\subset \grass(3, \wedge^2 V)$ be the
open locus of $2$-planes in $\proj(\wedge^2 V)$ 
which intersect $\grass(2,V)$ transversely in a nonsingular
conic curve.
\begin{lm}
\label{trick}
There is a canonical isomorphism $Z \eqq \cal{X}_{ruled}$.
\end{lm}
\bpf
The family of lines determined by 
a nonsingular 
plane conic $C \subset \grass(2,V) \subset \proj(\wedge^2 V)$
sweeps out an irreducible degree $2$ surface in $\proj(V)$.
There are three possibilities for this degree $2$ surface:
a double plane, a quadric cone, or a nonsingular
quadric surface. If the conic $C$ sweeps out a a double
plane $H\subset \proj(V)$, 
then $C \subset P \subset \grass(2,V)$ where
$P$ is the plane of all lines contained in $H$.
If $C$ sweeps out a quadric cone, then $C\subset P \subset
\grass(2,V)$ where $P$ is the plane of all lines
passing through the vertex of the cone. Hence, if
$C$ is the transverse intersection $P\cap \grass(2,V)$
of a plane, then $C$ must correspond to a ruling
of a unique nonsingular quadric surface. 
Conversely, a ruling
of a nonsingular quadric surface yields a conic curve
in $\grass(2,V)$ which is the transverse intersection of
a unique $2$-plane in $\proj(V)$. These maps
are easily seen to be algebraic.
\epf

By Lemma 1,
the ruled quadric $((Q),r)$ corresponds to 
a $3$-dimensional subspace $F\subset \wedge^2 V$.
$F$ is $\SO(V)$-invariant since $((Q),r)$ is stabilized
by $\SO(V)$. $F$ is therefore a $3$-dimensional representation
of $\SO(V)$. Let $s$ be the {\em other} ruling of $(Q)$.
$((Q),s)$ similarly corresponds to an invariant 
 $3$-dimensional subspace
$\tilde{F} \subset \wedge^2 V$. The $\SO(V)$ representation
$\wedge^2 V$ decomposes as $\wedge^2 V \eqq F \oplus \tilde{F}$.

\section{ $B\SO(V)$}
\label{bsofour}
Let $V\eqq \com^4$ be equipped with a non-degenerate
quadratic form as before. Approximations to
$E\SO(V)$ and $B\SO(V)$ are obtained via direct sums of the 
representation $V^*$.
Let $m>>0$ and let
$$W_m \subset \sumo_{1}^{m} V^*$$
denote the spanning locus. $W_m$ is the
locus of $m$-tuples of vectors of $ V^*$
which span $V^*$.
The natural action of $\SO (V)$ on
$W_m$ is free and has a geometric quotient
(see [P2]).
The codimension of the complement of
$W_m$ in $\sumo_{1}^{m}  V^*$ is $m-3$.
$W_m$ is an approximation of 
$E\SO (V)$ up to codimension $m-4$.
Therefore
$$B\SO (V)=\ \stackrel{Lim}{m \rarr \infty} \ W_m/\SO (V),$$
$$A^*_{\SO (V)}(\text{point})=\  \stackrel{Lim}{m \rarr \infty} \ 
A^*(W_m/\SO (V)).$$

There is a scalar $\com^*$-action on $W_m$.
Let $\proj(W_m)=W_m/\com^*$. Since this $\com^*$-action
commutes with the $\SO(V)$-action, there is diagram of
quotients:
\begin{equation}
\label{diax}
\begin{CD}
 W_m @>{\tau_1}>>  W_m/ \SO(V) \\
@V{i_1}VV @V{i_2}VV \\
\proj(W_m)= W_m/\com^* @>{\tau_2}>>\proj( W_m)/ \SO (V)\\
\end{CD}
\end{equation}
All the maps in (\ref{diax}) are quotient maps:
\begin{enumerate}
\item[(i)] $i_1$ is a free $\com^*$-quotient.
\item[(ii)] $i_2$ is a free $\com^*/(\pm)$-quotient.
\item[(iii)] $\tau_1$ is a free $\SO(V)$-quotient.
\item[(iv)] $\tau_2$ is a free $\proj\SO(V)$-quotient.
\end{enumerate} 
The existence of these quotients is easily deduced (see [P2]).

First consider the space $\proj(W_m)/ \SO(V)$.
Let $Q$ be the quadratic form on $V$ and let
$r$ be the ruling of the quadric surface $(Q)\subset  \proj(V)$
fixed in section \ref{rrr}.
An element $f \in \proj(W_m)$ yields a canonical embedding
$$\mu_f: \proj(V) \hookrightarrow \proj(\com^m).$$
The image under $\mu_f$ of $((Q),r)$ is 
a ruled quadric surface in $\proj(\com^m)$
associated canonically to $f \in \proj(W_m)$.
Since $\proj\SO(V)\subset \PGL(V)$ 
is exactly the stabilizer of the ruled quadric $((Q),r)$,
it follows that $\proj(W_m)/ \SO(V)$ is isomorphic
to the parameter space of ruled quadric surfaces in $\proj(\com^m)$.
Since a ruled quadric surface in $\proj(\com^m)$
 spans a unique $3$-plane  in
$\proj(\com^m)$, the parameter space is fibered over $\grass(4,m)$.
By Lemma \ref{trick}, the parameter space of ruled
quadric surface in $\proj(\com^m)$ is canonically 
isomorphic to an
open set $$Z \subset \grass(3, \wedge^2 S)$$
where $S\rarr \grass(4,m)$
is the tautological sub-bundle.

The Chow computations in section \ref{calcc} will require two 
results about line bundles. We have seen 
$\proj(W_m)/\SO(V)$ is canonically fibered over
$\grass(4,m)$. Let $c_1$ be the first Chern
class of the tautological bundle $S$ on $\grass(4,m)$.
Let $c_1$ also denote the pull-back of this class
to $\proj(W_m)/ \SO(V)$.
For $m>4$, $A^1(\proj(W_m)) \eqq \Z$ with generator
$c_1(\oh_{\proj}(-1))$ (which is the Chern class of the
line bundle associated
to the $\com^*$-bundle $i_1$.)
\begin{lm}
\label{lastl} $\tau_2^*(c_1)= c_1(\oh_{\proj}(-4))$.
\end{lm}
\bpf Elements of $\proj(W_m)$
correspond to embeddings of $\proj(V)$ in $\proj(\com^m)$.
The class $\tau_2^*(-c_1)$ is the divisor class
of embeddings that meet a fixed $(m-5)$-plane in $\proj(\com^m)=
\proj^{m-1}$. This divisor class is determined by a $4\times 4$
determinant. Hence, $\tau_2^*(-c_1)=c_1(\oh_{\proj}(4))$.
\epf

The map $i_2: W_m/ \SO(V) \rarr \proj(W_m)/ \SO(V)$
is a $\com^*/(\pm)$-bundle. 
Since there is an abstract isomorphism $\com^*/(\pm)\eqq \com^*$,
$i_2$ is also a $\com^*$-bundle. Let $N$ be the
line bundle on $\proj(W_m)/\SO(V)$ 
canonically associated to $i_2$.
\begin{lm}
\label{twot}
$\tau_2^*(N) \eqq \oh_{\proj}(-2)$.
\end{lm}
\bpf  
Let $i_1/(\pm): W_m/(\pm) \rarr \proj(W_m)$. The map
$i_1/(\pm)$ is a free $\com^*/(\pm)$-quotient.
The line bundle associated to the $\com^*/(\pm)$-bundle
$i_1/(\pm)$ is $\oh_{\proj}(-2)$.
The map $\tau_1/(\pm): W_m/(\pm) \rarr W_m/ \SO(V)$ is 
$\com^*/(\pm)$-equivariant. Hence, $\tau_2^*(N)\eqq
\oh_{\proj}(-2)$.
\epf

\section{\bf{ Chow Calculations}}
\label{calcc}
In this section,
the Chow ring of $W_m/ \SO(V)$ is determined
(up to codimension $m-4$). 
Consider the parameter space of ruled quadrics in
$\proj(\com^m)$:
$$Z \subset \grass(3,\wedge^2 S) \rarr \grass(4,m).$$
Let $D$ be the complement of $Z$ in $\grass(3, \wedge^2 S)$.
Following the notation of 
 section \ref{bsofour}, $W_m/\SO(V)$ is the
$\com^*$-bundle associated to a line bundle $N\rarr Z$. Therefore,
$$A^*(W_m/\SO(V)) \eqq A^*(Z)/ (c_1(N)) \eqq
A^*(\grass(3,\wedge^2 S))/(I_D, c_1(\barr{N}))$$
where $I_D \subset A^*(\grass(3,\wedge^2 S))$ is the ideal
generated by cycles supported on $D$ and
$\barr{N}$ is any extension of $N$ to $\grass(3,\wedge^2 S)$.
The ideal  $I_D$ is
determined by constructing a well-behaved variety which surjects
onto $D$. 

Let 
$\grass(2,S) \hookrightarrow \proj(\wedge^2 S)$
be the canonical relative Pl\"ucker embedding.
$D$ is exactly the locus of $2$-planes in the
the fibers of $\proj(\wedge^2 S)$ which do not
meet $\grass(2,S)$ transversely in a nonsingular
conic curve. Equivalently, $D$ is the locus
of $2$-planes $P$ in the fibers of $\proj(\wedge^2 S)$ which
satisfy one of the following conditions:
\begin{enumerate}
\item[(i)] $P\cap \grass(2,S)$ is a pair of distinct lines in $P$.
\item[(ii)] $P \cap \grass(2,S)$ is a double line in $P$.
\item[(iii)] $P \cap \grass(2,S)= P$.
\end{enumerate}

$D$ is dominated by a canonical Grassmannian bundle over
$\grass(2,S)$.
Let $B\rarr \grass(2,S)$ be the tautological sub-bundle.
By wedging, there is canonical surjective 
 bundle map on $\grass(2,S)$:
$$\wedge^2 S \otimes \wedge^2 B \rarr \wedge^4 S$$
which induces a 
canonical sequence on $\grass(2,S)$:
\begin{equation}
\label{aaa}
0 \rarr K \rarr \wedge^2 S \rarr \wedge^4 S 
\otimes (\wedge^2 B)^*\rarr
0.
\end{equation}
There is a canonical inclusion $\wedge^2 B \subset K$ and
a quotient sequence
\begin{equation}
\label{bbb}
0 \rarr \wedge^2 B \rarr K \rarr E \rarr 0
\end{equation}
on $\grass(2,S)$.
The geometric interpretation of these sequences is as
follows. Let $\xi\in \grass(2,S)$. $\proj(K_\xi)\subset\proj(
\wedge^2 S_\xi)$ is the projective tangent space to
$\grass(2, S_\xi)$ at $\xi$. $\proj(\wedge^2 B_\xi)$ in
$\proj(\wedge^2 S_\xi)$ is the Pl\"ucker image of the point
$\xi$. The fiber of the Grassmannian bundle
$$\grass(2,E) \rarr \grass(2,S)$$ over $\xi$
corresponds to the $2$-planes $P$ of $\proj(\wedge^2 S_{\xi})$
that are tangent to $\grass(2,S)$ at $\xi$.
There is a canonical map
$$\rho: \grass(2,E) \rarr D$$
which is a surjection of algebraic varieties.
Let $[P]\in D$. The fiber of $\rho$ over $[P]$
is simply the set of points of
$P \cap \grass(2,S)$ where $P$ is tangent to
$\grass(2,S)$. In case (i) above, the fiber is a point.
In case (ii), the fiber is a straight line in $\proj(\wedge^2 S)$.
In case (iii), the fiber is $2$-plane in $\proj(\wedge^2 S)$.
Hence, there is stratification of $D$ by intersection type (i-iii)
where $\rho$ is a projective bundle over each stratum.
The Chow groups of $\grass(2,E)$ therefore surject upon
the Chow groups of $D$.

The ideal $I_D$ is determined by calculating the
push-forwards of the Chow classes of $\grass(2,E)$ to
$\grass(3, \wedge^2 S)$. Consider the projection
$$\pi:\cal{G}= \grass(3, \wedge^2 S) \times
_{\grass(4,m)} \grass(2,S) \rarr \grass(3,\wedge^2 S).$$ 
The sequences (\ref{aaa}) and (\ref{bbb}) pull-back to
$\cal{G}$.
Let $F\rarr \grass(3, \wedge^2 S)$ denote the
tautological sub-bundle (and also let $F$ denote the
pull-back to $\cal{G}$ of this bundle).
There is a canonical inclusion
$$\iota: \grass(2,E) \hookrightarrow \cal{G}$$
determined by the sequences (\ref{aaa}) and (\ref{bbb}).
$\grass(2,E)\subset \cal{G}$ is the closed
subvariety of points $g \in \cal{G}$ where
$$\wedge^2 B_g \subset F_g \subset K_g.$$

The class $[\grass(2,E)]$ in Chow ring
of $A^*(\cal{G})$ is easily found by degeneracy calculations.
Let $c_1,c_2, c_3, c_4$ be the Chern classes of $S\rarr\grass(4,m)$.
Let $b_1, b_2$ be the Chern classes of $B \rarr \grass(2,S)$.
Let $f_1,f_2,f_3$ be the Chern classes of 
$F\rarr \grass(3,\wedge^2 S)$. Since $\cal{G}$ is
a tower of Grassmannian bundles, these Chern classes
$c_i, b_j, f_k$ generate $A^*(\cal{G})$.
Let $Y$ be 
the locus of points $g\in \cal{G}$ such that $F_g \subset K_g$.
$Y$
is the nonsingular degeneracy locus of the canonical
bundle map on $\cal{G}$,
$$ F  \rarr \wedge^4 S \otimes (\wedge^2 B)^*,$$
obtained from the inclusion $F \subset \wedge^2 S$ and
sequence (\ref{aaa}).
By the Thom-Porteous formula on $\cal{G}$ (see [F]),
\begin{eqnarray*}
A^*(\cal{G})\ni
[Y] &= &c_3(F^*\otimes \wedge^4 S \otimes (\wedge^2 B)^*) \\
 &=& -f_3 +(c_1-b_1)f_2- (c_1-b_1)^2 f_1+(c_1-b_1)^3.
\end{eqnarray*}
$Y$ is canonically isomorphic to 
the  Grassmannian bundle $\grass(3, K) \rarr \grass(2,S)$.
There is natural bundle quotient sequence on $Y$:
\begin{equation}
\label{ccc}
0 \rarr F \rarr K \rarr K/F \rarr 0.
\end{equation}
The locus $\grass(2,E)\subset Y$ is the set 
of points $y\in Y$ such that
$\wedge^2 B_y \subset F_y$. $\grass(2,E)\subset Y$
is the nonsingular degeneracy
locus of the canonical bundle map on $Y$,
$$\wedge^2 B \rarr K/F,$$
obtained from the sequences (\ref{bbb}) and (\ref{ccc}).
By the Thom-Porteous formula on $Y$,
\begin{eqnarray*}
A^*(Y)\ni 
[\grass(2,E)] & = &c_2( (K/F)\otimes (\wedge^2 B)^*) \\
& =& b_1^2-c_1b_1+ c_1^2-2c_1 f_1+ f_1^2-f_2+2c_2.
\end{eqnarray*}
The class $[\grass(2,E)] \in A^*(\cal{G})$ is there
expressed by 
\begin{eqnarray*}
A^*(\cal{G}) \ni [\grass(2,E)] & = &
\big(-f_3 +(c_1-b_1)f_2- (c_1-b_1)^2 f_1+(c_1-b_1)^3\big)  \\
& &
\cdot \big(b_1^2-c_1b_1+ c_1^2-2c_1 f_1+ f_1^2-f_2+2c_2\big).
\end{eqnarray*}

Since $\grass(2,E)$ is a Grassmannian bundle over $\grass(2,S)$,
the Chow ring of $\grass(2,E)$ is generated over $A^*(\grass(2,S))$
by the Chern classes $h_1, h_2$ of the tautological
sub-bundle $H\rarr \grass(2,E)$. Via
the embedding  $\iota:\grass(2,E) \hookrightarrow
\cal{G}$, $H$ is isomorphic to  $\iota^*(F)/ 
\iota^*(\wedge^2 B)$. 
The  Chern classes
$h_1$ and $h_2$ can be expressed via 
$\iota$ in terms of the
classes $b_j$ and $f_k$.
Therefore, the classes 
\begin{equation*} 
\grass(2,E) \cap M(c_1,c_2,c_3,c_4,b_1,b_2,f_1,f_2,f_3)
\end{equation*}
(where $M$ is monomial in the Chern classes) span the
Chow ring of $\grass(2,E)$. The ideal $I_D \subset A^*(\grass(3,
\wedge^2 S))$ is generated by the $\pi$
push-forwards of the classes (\ref{clazz}): 
\begin{equation}
\label{clazz} 
[\grass(2,E)] \cdot M(c_1,c_2,c_3,c_4,b_1,b_2,f_1,f_2,f_3) \in 
A^*(\cal{G})
\end{equation}
Since the classes $c_i, f_k$ in $A^*(\cal{G})$ are
pull-backs from $\grass(3,\wedge^2 S)$, $I_D$ is generated
by the elements
$$\pi_*([\grass(2,E)] \cdot M(b_1,b_2)).$$
By the standard relations satisfied by the classes
$b_1$ and $b_2$ over $A^*(\grass(4,m))$, it follows that
$I_D$ is generated by:
$$\pi_*([\grass(2,E)])$$
$$\pi_*([\grass(2,E)]\cdot b_1)$$
$$\pi_*([\grass(2,E)]\cdot b_1^2) \ \ \ \pi_*([\grass(2,E)]\cdot b_2)$$
$$\pi_*([\grass(2,E)]\cdot b_1b_2)$$
$$\pi_*([\grass(2,E)]\cdot b^2_1b_2).$$

Since the class $[\grass(2,E)]$ is determined explicitly in
$A^*(\cal{G})$,
these six push-forwards can be easily computed by hand or by the
the Maple package Schubert ([KS]).
The first push-forward is:
 $$\pi_*([\grass(2,E)])=13c_1-2f_1.$$
\begin{lm}
\label{killl}
The pair $(13c_1-2f_1, c_1(\barr{N}))$ generates
$A^1(\grass(3,\wedge^2 S))$.
\end{lm}
\bpf
Recall the notation of diagram (\ref{diax}):
$$\tau_2: \proj(W_m) \rarr Z \subset \grass(3,\wedge^2 S).$$
Let $L=c_1(\oh_{\proj}(-1))$ be a generator of
$A^1(\proj(W_m))$.
By Lemma 2, $\tau_2^*(c_1)\eqq 4L$.
Since $\tau_2^*([D])=\tau_2^*(13c_1-2f_1)=0$,
$\tau_2^*(f_1)=26 L$. Therefore the
image of $\tau_2^*$ is the subgroup $\Z(2L)$.

Since $[D]=13c_1-2f_1$ is not divisible in 
$A^1(\grass(3,\wedge^2 S))$, $A^1(Z)\eqq \Z$ and
$\tau_2^*$ is an isomorphism:
$$\tau_2^*: A^1(Z) \stackrel{\sim}{\rarr} \Z(2L).$$
By Lemma 3, $\tau_2^*(c_1(N))=2L$. Therefore,
$c_1(N)$ generates $A^1(Z)$. It now follows that the pair
$(13c_1-2f_1, c_1(\barr{N}))$ generates the group
$A^1(\grass(3,\wedge^2 S))$.
\epf

\noindent
Therefore, $(I_D, c_1(\barr{N}))=(I_D,c_1, f_1).$

By Lemma \ref{killl}, it suffices to compute 
the five remaining push-forwards
modulo the ideal $J=(c_1,f_1)$. The 
results are (modulo $J$):
\begin{enumerate}
\item[{}] $\pi_*([\grass(2,E)]\cdot b_1)=0.$
\item[{}] $\pi_*([\grass(2,E)]\cdot b_1^2)=-2f_3.$
\item[{}] $\pi_*([\grass(2,E)]\cdot b_2)=c_3-f_3.$
\item[{}] $\pi_*([\grass(2,E)]\cdot b_1b_2)=(c_2-f_2)^2-4c_4.$
\item[{}] $\pi_*([\grass(2,E)]\cdot b_1^2 b_2)=c_2f_3+f_2c_3.$
\end{enumerate}
Hence $(I_D, c_1(\barr{N}))=(c_1,f_1,2c_3, c_3-f_3,
(c_2-f_2)^2-4c_4, (c_2-f_2)c_3)$.

The ring $A^*(\grass(4,m))$ is freely generated
(up to codimension $m-4$) by $c_1,c_2,c_3,c_4$.
The ring $A^*(\grass(3,\wedge^2 S))$ has the
following presentation (up to codimension $m-4$):
$$A^*(\grass(3,\wedge^2 S))\eqq \Z[c_1,c_2,c_3,c_4, f_1,f_2,f_3]/
(t_4,t_5,t_6)$$
where the $t_4,t_5, t_6$ are the Chern classes of the tautological
quotient bundle $T$: $0 \rarr F \rarr \wedge^2 S  \rarr T \rarr 0.$
We find (modulo $J$):
\begin{enumerate}
\item[{}] $t_4= (c_2-f_2)^2-4c_4$.
\item[{}] $t_5= 2f_2f_3-2c_2f_3$.
\item[{}] $t_6= f_2(-(c_2-f_2)^2+4c4)+f_3^2-c_3^2$.
\end{enumerate}
There is a presentation:
$A^*(\grass(3,\wedge^2 S))/(I_D,c_1(\barr{N})) \eqq$
\begin{equation}\label{prezz} 
\Z[c_1,c_2,c_3,c_4,f_1,f_2,f_3]/(I_D, c_1, f_1, t_4,t_5,t_6)
\end{equation}
(up to codimension $m-4$).
Surprisingly, the relations $t_4,t_5,t_6$
are contained in the ideal $(I_D, c_1,f_1)$. By the limit
procedure,
$$A_{\SO(4)}^*(\text{point})\eqq
\Z[c_1,c_2,c_3,c_4, f_2]/(c_1,2c_3,(c_2-f_2)c_3, (c_2-f_2)^2-4c_4).$$
The vector bundles $S$, $F \subset \wedge^2 S$ on the
approximation $W_m/\SO(V)$ are easily seen to be obtained
from the principal $\SO(V)$-bundle
$$W_m \rarr W_m/ \SO(V)$$
and the representations $V$, $F \subset \wedge^2 V$ defined
in section \ref{rrr}. Define $$x=c_2-f_2.$$ Theorem 1 is
proved.

\vspace{+10 pt}
\noindent
Department of Mathematics \\
University of Chicago \\
5734 S. University Ave. 60637 \\
rahul@math.uchicago.edu

\begin{thebibliography}{[EG2]}
\bibitem[EG1]{a} D. Edidin and W. Graham, {\em Equivariant intersection
                theory}, preprint 1996.
\bibitem[EG2]{aa} D. Edidin and W. Graham, {\em Characteristic
classes of principal bundles in algebraic geometry}, preprint 1995.
\bibitem[F]{b} W. Fulton, {\em Intersection theory}, 
             Springer-Verlag:
             Berlin, 1984.
\bibitem[FH]{bb} W. Fulton and J. Harris, {\em Representation
Theory}, Springer-Verlag: New York, 1991.
\bibitem[KS]{bbb} S. Katz and S. A. Stromme,
{\em Schubert: A Maple package for intersection theory in
algebraic geometry}.
\bibitem[P1]{d} R. Pandharipande, {\em The Chow ring of the
 non-linear Grassmannian}, preprint 1996.
\bibitem[P2]{e} R. Pandharipande, {\em The Chow ring of the
Hilbert scheme of rational normal curves}, preprint 1996.
\bibitem[T] {f} B. Totaro, {\em The Chow ring of the
symmetric group}, preprint 1994.
\end{thebibliography}
\end{document}